\documentclass[a4paper]{amsart}
\usepackage{amsthm,amsmath,amssymb}
\usepackage[latin1]{inputenc}
 \newtheorem{thm}{Theorem}
 \newtheorem{prop}[thm]{Proposition}
 \newtheorem{lemma}[thm]{Lemma}
 \newtheorem{kor}[thm]{Corollary}
 
 \theoremstyle{definition}

 \theoremstyle{remark}
 \newtheorem{remark}[thm]{Remark}

 \def\Spec{{\rm Spec}}

\def\Mat{{\rm Mat}}

 \def\id{{\rm id}}

 \def\sb{{\rm sb}}
 \def\ad{{\rm ad}}
 \def\det{{\rm det}}

 \def\dom{{\rm dom}}
 \def\M{{M}}
\def\a{{a}}
 \def\P{{P}}
 \def\N{{N}}
 
\begin{document}
\begin{title}
{Connected components of closed affine Deligne-Lusztig varieties} 
\end{title}
\author{Eva Viehmann}
\address{Mathematisches Institut der Universit\"{a}t Bonn\\ Beringstrasse 1\\53115 Bonn\\Germany}
\date{}
\begin{abstract}{We determine the set of connected components of closed affine Deligne-Lusztig varieties for special maximal compact subgroups of split connected reductive groups. We show that there is a transitive group action on this set. Thus such an affine Deligne-Lusztig variety has isolated points if and only if its dimension is $0$. We also obtain a description of the set of these varieties that are zero-dimensional.}
\end{abstract}
\maketitle
\section{Introduction}\label{secdlvcontext}
Let $k$ be a finite field with $q=p^r$ elements and let $\overline{k}$ be an algebraic closure. Let $F=k((t))$ and let $L=\overline{k}((t))$. Let $\mathcal{O}_F$ and $\mathcal{O}_L$ be the valuation rings. We denote by $\sigma:x\mapsto x^q$ the Frobenius of $\overline{k}$ over $k$ and also of $L$ over $F$.

Let $G$ be a split connected reductive group over $\mathcal{O}_F$ and let $A$ be a split maximal torus. Let $B$ be a Borel subgroup containing $A$. For dominant elements $\mu,\mu'\in X_*(A)$ we say that $\mu'\preceq \mu$ if $\mu-\mu'$ is a non-negative linear combination of positive coroots. We write $\mu_{\dom}$ for the dominant element in the orbit of $\mu\in X_*(A)$ under the Weyl group $W$ of $A$ in $G$. For $\mu\in X_*(A)$ and $x\in L^{\times}$ we denote by $x^{\mu} \in A(L)$ the image of $x$ under the homomorphism $\mu: \mathbb{G}_m\rightarrow A$.

We recall the definitions of affine Deligne-Lusztig varieties and closed affine Deligne-Lusztig varieties from \cite{Rapoport1}, \cite{GHKR}. Let $K=G(\mathcal{O}_L)$ and let $X=G(L)/K$ be the affine Grassmannian. For $b\in G(L)$ and a dominant coweight $\mu\in X_*(A)$ the affine Deligne-Lusztig variety $X^G_{\mu}(b)=X_{\mu}(b)$ is the locally closed reduced $\overline{k}$-subscheme of $X$ defined by 
\begin{equation*}
X_{\mu}(b)(\overline{k})=\{g\in G(L)/K\mid g^{-1}b\sigma(g)\in Kt^{\mu}K\}.
\end{equation*} 
The closed affine Deligne-Lusztig variety is the closed reduced subscheme of $X$ defined by 
\begin{equation*}
X_{\preceq \mu}(b)=\bigcup_{\mu'\preceq \mu}X_{\mu'}(b).
\end{equation*} Both $X_{\mu}(b)$ and $X_{\preceq\mu}(b)$ are locally of finite type.

By $[x]$ we denote the $\sigma$-conjugacy class of an element $x\in G(L)$. Left multiplication by $g\in G(L)$ induces an isomorphism between $X_{\mu}(b)$ and $X_{\mu}(gb\sigma(g)^{-1})$. Thus the isomorphism class of the affine Deligne-Lusztig variety only depends on $[b]$ and not on $b$.

We write $\pi_1(G)$ for the quotient of $X_*(A)$ by the coroot lattice of $G$. In \cite{Kottwitz2}, Kottwitz defines a homomorphism $\kappa_G: G(L)\rightarrow \pi_1(G)$ which induces a locally constant map $\kappa_G:X\rightarrow \pi_1(G)$. For $G=GL_h$ we have $\pi_1(G)\cong\mathbb{Z}$ and $\kappa_G(g)=v_t(\det~g)$ where $v_t$ denotes the usual valuation on $L=\overline{k}((t))$.

Let $\nu$ be the Newton vector associated to $b$, compare \cite{Kottwitz1}. It is a dominant element of $X_*(A)_{\mathbb{Q}}$. In \cite{KottRapo} Kottwitz and Rapoport give a criterion for $X_{\mu}(b)$ to be nonempty, compare also \cite{GHKR}, Proposition 5.6.1. From now on we assume that this is the case. Then $\kappa_G(b)=\mu$ and $\mu-\nu$ is a positive linear combination of positive coroots with rational coefficients.

Let $\P$ be a standard parabolic subgroup of $G$. Then $\P=\M \N$, where $\N$ is the unipotent radical of $\P$ and where $\M$ is the unique Levi subgroup of $\P$ containing $A$. Applying the construction of $\kappa$ to $\M$ rather than $G$ we obtain a homomorphism $\kappa_\M:\M(L)\rightarrow \pi_1(\M)$. The inclusion $\M(L)/\M(\mathcal{O}_L)\hookrightarrow G(L)/G(\mathcal{O}_L)$ induces for each $\mu$ and each $b\in \M(L)$ an inclusion $X_{\mu}^{\M}(b)\hookrightarrow X_{\mu}^G(b)$. Here $X_{\mu}^\M(b)$ denotes the affine Deligne-Lusztig variety for $\M$. 

Let $A_{\P}$ denote the identity component of the center of $\M$. Let
\begin{equation*}
\mathfrak{a}_{\P}^+=\{x\in X_*(A_{\P})\otimes_{\mathbb{Z}}\mathbb{R}\mid \langle\alpha,x\rangle >0 \text{ for every root }\alpha \text{ of }A_{\P} \text{ in }\N\}.
\end{equation*} 
As in \cite{Kottwitz2} let $\P_b=\M_b\N_b$ be the unique standard parabolic subgroup of $G$ such that $M_b$ is the centralizer of $\nu$. Then the $\sigma$-conjugacy class $[b]$ contains an element $b'$ with the following properties: $b'$ is basic in $\M_b$ and $\kappa_{\M_b}(b')$, considered as an element of $X_*(A_{\P_b})\otimes_{\mathbb{Z}}\mathbb{R}$, lies in $\mathfrak{a}_{\P_b}^+$. We assume that $b=b'$. 

The proof of the Hodge-Newton decomposition by Kottwitz (see \cite{Kottwitz3}) yields the following result. 
\begin{thm}\label{thmkottwitz}
Let $\P=\M \N\subseteq G$ be a standard parabolic subgroup with $\P_b\subseteq \P$. If $\kappa_{\M}(b)=\mu$, then the morphism $X_{\mu}^{\M}(b)\hookrightarrow X_{\mu}^G(b)$ is an isomorphism. 
\end{thm}
Indeed, the proof given by Kottwitz can be taken almost literally, once replacing his condition by ours, to show this stronger theorem. For more details, especially for the corresponding statement for unramified groups $G$, compare \cite{hn}.

We call a pair $(\mu,b)$ indecomposable with respect to the Hodge-Newton decomposition if for all standard parabolic subgroups $\P$ with $\P_b\subseteq \P=\M\N\subsetneq G$ we have $\kappa_{\M}(b)\neq \mu$. Given $G$, $\mu$, and $b$, we may always pass to a Levi subgroup $\M$ of $G$ in which $(\mu,b)$ is indecomposable. For a description of the affine Deligne-Lusztig varieties it is therefore sufficient to consider pairs $(\mu,b)$ which are indecomposable with respect to the Hodge-Newton decomposition. 

Let $G_{\ad}$ be the adjoint group of $G$. We denote the images of $b$ and $\mu$ in $G_{\ad}$ also by $b$ and $\mu$. Then the sets of connected components of $X^G_{\preceq\mu}(b)$ and $X^{G_{\ad}}_{\preceq\mu}(b)$ are closely related. More precisely, we prove in Lemma \ref{lemnormgm} that $$\pi_0(X^{G_{\ad}}_{\preceq\mu}(b))\cong \pi_0(X^{G}_{\preceq\mu}(b))/Z(G)(F).$$ Here $Z(G)$ is the center of $G$. Using this we show that it is enough to describe the set of connected components in the case that $G$ is simple.

Let $$J=\{g\in G(L)\mid g\circ b\sigma=b\sigma\circ g\}.$$ Then there is a canonical $J$-action on $X_{\mu}(b)$ for each $\mu$. We prove that already a subgroup of $J$ acts transitively on the set of connected components of $X_{\preceq\mu}(b)$. From this we obtain the following theorem.

\begin{thm}\label{thmzshk}
Let $G$, $\mu$, and $b$ be as above and indecomposable with respect to the Hodge-Newton decomposition. Assume that $G$ is simple.
\begin{enumerate}
\item Either $\kappa_{\M}(b)\neq \mu$ for all proper standard parabolic subgroups $\P$ of $G$ with $b\in \M$ or $[b]=[t^{\mu}]$ with $t^{\mu}$ central. 
\item In the first case, $\kappa_{G}$ induces a bijection $\pi_0(X_{\preceq\mu}(b))\cong\pi_1(G)$. 
\item In the second case, $X_{\mu}(b)=X_{\preceq \mu}(b)\cong J/(J\cap K)\cong G(F)/G(\mathcal{O}_F)$ is discrete.
\end{enumerate} 
\end{thm}

The transitivity of the $J$-action on the set of connected components implies that the map $\pi_0(X_{\mu}(b))\rightarrow \pi_0(X_{\preceq\mu}(b))$ induced by the inclusion is surjective. In Section \ref{secopen} we give an example to show that in general the map is not injective. There are also examples where the action of $J$ on $\pi_0(X_{\mu}(b))$ is not transitive.

Using the dimension formula for affine Deligne-Lusztig varieties from \cite{GHKR}, \cite{dimdlv} we show the following theorem. This characterisation of zero-dimensional affine Deligne-Lusztig varieties has been conjectured by Rapoport.
\begin{thm}\label{thmdim0}
Let $(\mu,b)$ be indecomposable with respect to the Hodge-Newton decomposition. Then the following assertions are equivalent.
\begin{enumerate}
\item $X_{\mu}(b)$ has an isolated point
\item $\dim X_{\mu}(b)= 0$
\item Let $G_1\times\dotsm\times G_n$ be the decomposition of $G_{\ad}$ into simple factors. Then for each pair $(\mu_i,b_i)$ of images of $\mu$ and $b$ in some $G_i$, we have either that the $\sigma$-conjugacy class $[b_i]=[t^{\mu_i}]$ is central in $G_i$ or that $G_i=PGL_h$ for some $h$, that $b_i$ is basic and that $\mu_i\in X_*(A_i)\cong \mathbb{Z}^h/\mathbb{Z}$ is minuscule of the form $(0,\dotsc,0,1)$ or $(0,1,\dotsc,1)$. Here, the torus $A_i\subseteq G_i$ is the image of $A$ in $G_i$.
\end{enumerate}
\end{thm}

For $G=GL_n$, Theorems \ref{thmzshk} and \ref{thmdim0} are part of my thesis, compare \cite{diss}, 2. The proof given there is along the same lines, but more explicit.\\

\noindent{\it Acknowledgement.} I thank M. Rapoport for introducing me to this subject and for his interest in my work. I am grateful to R. Kottwitz for valuable discussions. This paper was written during a stay at the Universit\'{e} Paris-Sud at Orsay which was supported by a fellowship within the Post-Doc program of the German Academic Exchange Service (DAAD). I thank the Universit\'{e} Paris-Sud for its hospitality.

\section{Connected components}\label{secconncomp}
\subsection{Notation and preliminary reductions}\label{secnot}
We recall briefly some of the notation associated with Kottwitz' classification of isocrystals, see \cite{Kottwitz2}. For split groups,  $\kappa=\kappa_G:G(L)\rightarrow \pi_1(G)$ has the following easy description. Let $b\in G(L)$ and $r_B(b)=\mu\in X_*(A)$ the unique element such that $b\in Kt^{\mu}U(L)$ where $U$ is the unipotent radical of $B$. Then $\kappa(b)$ is the image of $\mu$ under the canonical projection from $X_*(A)$ to $\pi_1(G)$. 

Recall that $M_b\subseteq G$ is the centralizer of $\nu$ and that we assume that $b\in M_b$. Then $M_b$ is the largest standard Levi subgroup of $G$ in which $b$ is basic. The $\sigma$-conjugacy class of a basic element is determined by its value of $\kappa$ (see \cite{Kottwitz2}). Especially, the $\sigma$-conjugacy class $[b]$ of $b$ is determined by $\nu$ and $\kappa_{M_b}(b)$. Note that the quotient of the coroot lattices of $G$ and $M_{b}$ is generated by the simple coroots of $A$ in $N_{b}$, and hence torsion free. The images of $\nu$ and $\kappa_{M_b}(b)$ in $\pi_1(M_{b})\otimes\mathbb{Q}$ are equal. Thus $\nu$ and $\kappa_G(b)$ determine $\kappa_{M_{b}}(b)$ and also $[b]$.

In \cite{GHKR}, 5.9 it is shown that there exists a standard parabolic subgroup $P_{\sb}$ of $G$ with Levi factor $M_{\sb}$ containing $A$, unipotent radical $N_{\sb}$ and the following properties. The intersection $[b]\cap M_{\sb}(L)$ is nonempty, thus we may assume $b\in M_{\sb}(L)$. Furthermore, $b$ is superbasic in $M_{\sb}$, i.e. no $\sigma$-conjugate of $b$ lies in a proper Levi subgroup of $M_{\sb}$. The root system of $M_{\sb}$ is then a disjoint union of root systems of type $A_l$. 

In the following we construct a special representative of $[b]$. We have $M_{\sb,\ad}= \prod_i M_i$ with $M_i\cong PGL_{h_i}$ for some $h_i$. We may assume that the maximal torus $A_i$ and the Borel subroup $B_i$ of $M_i$ that are induced by $A$ and $B$ are the subgroups of diagonal and of upper triangular matrices in $PGL_{h_i}$. Let $W_{M_i}$ be the Weyl group of $A_i$ in $M_i$. Let $(b_i)$ be the image of $b$ in $\prod_i M_i$. Let $e_j^i$ with $j=0,\dotsc,h_i-1$ be the standard basis of $L^{h_i}$. For $j\in\mathbb{Z}$ define $e_j^i$ by $e_{j+h_i}^i=te_j^i$. Then $[b_i]$ contains a representative with $b_i(e_j^i)=e_{j+m_i}^i$ for all $j$ and some $m_i\in\mathbb{Z}$ with $(m_i,h_i)=1$. Especially $b_i=t^{\mu_{M_i}}w_i$ with $\mu_{M_i}\in X_*(A_i)$ minuscule and dominant and $w_i\in W_{M_i}$. Here $w_i(e_j^i)$ for $j\in [0,h_i-1]$ is equal to  $e_{j+m_i-lh_i}^i$ where $l$ is such that the index is again in $[0,h_i-1]$. The second step is to use these special $b_i$ to define a similar normal form for $b$. As $(\mu_{M_i})_i=\kappa_{M_{\sb}}(b)$ as elements of $\pi_1(M_{\sb})/X_*(A_{P_{\sb}})$, there is some $\mu_{M_{\sb}}\in X_*(A)$ mapping to $(\mu_{M_i})_i$ in $X_*(A/A_{P_{\sb}})$ and with $\kappa_{M_{\sb}}(b)=\mu_{M_{\sb}}$. Let $w=(w_i)\in \prod W_{M_i}\subseteq W$. Then $\kappa_{M_{\sb}}(t^{\mu_{M_{\sb}}}w)=\kappa_{M_{\sb}}(b)$. As $t^{\mu_{M_{\sb}}}w$ and $b$ are basic in $M_{\sb}$, they are $\sigma$-conjugate. Hence we may assume that $b=t^{\mu_{M_{\sb}}}w$. In Corollary \ref{corsurj} we show that $(\mu_{M_{\sb}})_{\dom}=\mu_{\min}$ is minimal among the dominant elements $\mu$ with $X_{\mu}(b)\neq \emptyset$.

From the special form of $b$ we obtain that $\id\in X^{M_{\sb}}_{\mu_{M_{\sb}}}(b)$. Thus it is in $X^G_{\mu_{\min}}(b)$. As $J$ commutes with $b\sigma$, its image $J/(K\cap J)$ in $X$ is also contained in $X^G_{\mu_{\min}}(b)$.

Recall that $A_G$ denotes the identity component of the center $Z(G)$ of $G$.
\begin{lemma}\label{lemnormgm}
For $n\in\pi_1(G)\cong \pi_0(X)$ let $X^G_{\preceq\mu}(b)_n$ be the intersection of $X_{\preceq\mu}(b)$ and the connected component of $X$ corresponding to $n$. Then the $X^G_{\preceq\mu}(b)_n$ for the different $n$ are isomorphic. 

Let $G_{\ad}$ be the adjoint group of $G$ and denote the images of $b$ and $\mu$ in $G_{\ad}$ again by $b$ and $\mu$. Let $X_{\preceq\mu}^{G_{\ad}}(b)_0$ be the intersection of $X_{\preceq\mu}^{G_{\ad}}(b)$ and the connected component of the affine Grassmannian of $G_{\ad}$ associated to $0\in\pi_1(G_{\ad})$. Recall that we assume that $X^G_{\preceq\mu}(b)\neq\emptyset$. Then $X_{\preceq\mu}^{G}(b)_0\cong X_{\preceq\mu}^{G_{\ad}}(b)_0$.
\end{lemma}
\begin{proof}
The map $\kappa_{M_{\sb}}|_{J\cap M_{\sb}}:J\cap M_{\sb}\rightarrow \pi_1(M_{\sb})$ is surjective. Indeed, to show this, we may assume that $M_{\sb}$ is adjoint (as $Z(M_{\sb})(F)\subseteq J\cap M_{\sb}$) and consider each simple factor separately. Such a factor $M_i\cong PGL_{h_i}$ has fundamental group $\mathbb{Z}/h_i\mathbb{Z}$. An element of $J\cap M_{\sb}$ mapping to $n\in \mathbb{Z}/h_i\mathbb{Z}$ is then explicitly given by $e_j^i\mapsto e_{j+n}^i$ for all $i$. As $\pi_1(G)$ is a quotient of $\pi_1(M_{\sb})$, this implies that $J$ is mapped surjectively to $\pi_1(G)$. Especially, it defines ismorphisms between the different $X^G_{\preceq\mu}(b)_n$ for $n\in \pi_1(G)$.

For the second assertion let $g\in X^{G_{\ad}}_{\preceq\mu}(b)_0$ and $\tilde{g}\in G(L)$ in the inverse image of $g$ under the surjection $G(L)\cap \ker(\kappa_G)\rightarrow G_{\ad}(L)\cap \ker(\kappa_{G_{\ad}})$. We have $Z(G)(L)\subseteq K\{t^{\alpha}\mid\alpha\in X_*(A_G)\}$. Thus $\tilde{g}^{-1}b\sigma(\tilde{g})\in Kt^{\mu}K\cdot\{t^{\alpha}\mid\alpha\in X_*(A_G)\}$. But as $X^G_{\preceq\mu}(b)\neq \emptyset$, a set $X^G_{\preceq\mu+\alpha}(b)$ for some central $\alpha$ can only be nonempty if $\mu=\kappa(b)=\mu+\alpha$ in $\pi_1(G)$, hence if $\alpha=0$. Thus $\tilde{g}\in X^G_{\preceq\mu}(b)$. This implies that the map is surjective. Let now $\tilde{g}'$ be a second inverse image of $g$. Then $\tilde{g}'\in \tilde{g}t^{\alpha}K$ for some $\alpha\in X_*(A_G)$. Hence the map $X^G_{\mu}(b)_0\rightarrow X^{G_{\ad}}_{\mu}(b)_0$ is a bijection. 
\end{proof}
\begin{remark}\label{rem0}
To determine $\pi_0(X^{G}_{\preceq\mu}(b))$ for a split connected reductive group $G$, we will determine the fibers of the map $\kappa:\pi_0(X^{G}_{\preceq\mu}(b))\rightarrow \pi_1(G)$. Note that $\kappa$ is induced by a homomorphism $G(L)\rightarrow \pi_1(G)$ and hence it is compatible with the $J$-action. By the preceding lemma we may thus assume that $G$ is adjoint. Then the closed affine Deligne-Lusztig variety is a product of the corresponding varieties associated to the simple factors of $G$. Hence for the computation of the set of connected components it is no restriction to assume that $G$ is simple.
\end{remark}

\subsection{Proof of Theorem \ref{thmzshk} (1)}
As $\kappa_G(b)=\mu$, the $\sigma$-conjugacy classes of $b$ and $t^{\mu}$ are equal if and only if $\nu=\mu$. If this is the case, the indecomposability applied to $M_b$ implies that $M_b=G$. Thus $\langle\alpha,\mu\rangle=0$ for every positive root $\alpha$ of $G$. Hence $t^{\mu}$ is central. It remains to prove $\nu=\mu\in X_*(A)_{\mathbb{Q}}$ using the assumption that there is a proper standard parabolic subgroup $P=MN$ of $G$ with $\kappa_M(b)=\mu$.

We use induction on the distance in the Dynkin diagram of $G$ between the simple root $\alpha$ and a simple root of $A$ in $N$ to show the following assertion.

{\bf Claim.} {\it If we write $\mu-\nu$ as a linear combination of the simple coroots of $G$, the coefficients of the coroot associated to $\alpha$ and all its neighbours in the Dynkin diagram vanish.}

In particular, $\langle\alpha,\mu-\nu\rangle=0$. If this claim is true for all $\alpha$, we have $\mu=\nu$, and Theorem \ref{thmzshk} (1) follows.

We first show that $\langle \alpha,\nu\rangle=0$ for all simple roots $\alpha$ of $A$ in $N$. As $\nu$ is dominant, all these products are non-negative. Assume that $\langle \alpha,\nu\rangle>0$ for some $\alpha$. Let $P'=M'N'\supseteq P$ be the maximal standard parabolic subgroup of $G$ corresponding to $\alpha$. Then $b\in M_b(L)\subseteq P'(L)$. From $\kappa_M(b)=\mu$ and $M\subseteq M'$ we obtain $\kappa_{M'}(b)=\mu$. But as $M_b\subseteq M'$, this is a contradiction to the indecomposability with respect to the Hodge-Newton decomposition.

As $\mu$ is dominant,
\begin{equation}\label{glmunu}
\langle\alpha,\mu-\nu\rangle= \langle\alpha,\mu\rangle+0\geq 0.
\end{equation}
The assumption $\kappa_M(b)=\mu$ implies that $\mu-\nu$ is a linear combination of simple coroots of $M$ with nonnegative rational coefficients. Hence if we consider it as a linear combination of the simple coroots of $G$, the coefficient of $\alpha^{\vee}$ vanishes. Thus the left hand side of (\ref{glmunu}) is a non-negative linear combination of the (non-positive) products of $\alpha$ with the other simple coroots. Hence, all coefficients in $\mu-\nu$ corresponding to neighbours of $\alpha$ in the Dynkin diagram vanish.

For the induction step assume that the assertion is true for all simple roots of $G$ whose distance to a simple root of $A$ in $N$ is at most $d_0$. Let $\alpha$ be a simple root of $G$ whose distance to the simple roots in $N$ is $d_0+1$. Thus the induction hypothesis implies that the coefficient of $\alpha^{\vee}$ in $\mu-\nu$ vanishes. If $\langle\alpha,\nu\rangle=0$, the same argument as above shows that the claim also holds for $\alpha$. Assume $\langle\alpha,\nu\rangle>0$. Let $P'=M'N'$ be the maximal standard parabolic subgroup of $G$ corresponding to $\alpha$. Then $\langle\alpha,\nu\rangle>0$ implies that $b\in M_b\subseteq M'$. Consider $r_{B\cap M'}^{M'}(b)-\nu$, which is a rational linear combination of coroots of $M'$. Here $r_{B\cap M'}^{M'}$ is as in the definition of $\kappa$ in Section \ref{secnot}. The induction hypothesis implies that $\mu-\nu$ is also a rational linear combination of coroots of $M'$. Hence $r_{B\cap M'}^{M'}(b)-\mu\in X_*(A)_{\mathbb{Q}}$ is a rational linear combination of coroots of $M'$. But as $\mu=\kappa_G(b)=r_{B\cap M'}^{M'}(b)$ in $\pi_1(G)$, the difference $\mu-r_{B\cap M'}^{M'}(b)\in X_*(A)$ is an integral linear combination of coroots of $G$. Using that the quotient of the two coroot lattices is torsion free, we obtain that $\mu-r_{B\cap M'}^{M'}(b)\in X_*(A)$ is an integral linear combination of coroots of $M'$. Therefore $\kappa_{M'}(b)=\mu$. Together with $M_b\subseteq M'$, this is a contradiction to the irreducibility with respect to the Hodge-Newton decomposition.\qed

\subsection{Proof of Theorem \ref{thmzshk} (3)}
Note that the map $k\mapsto k^{-1}\sigma(k)$ from $K$ to $K$ is surjective.  
Using this and that $[b]=[t^{\mu}]$ is central we obtain
\begin{eqnarray*}
X_{\mu}(b)&\cong&\{g\in G(L)/K\mid g^{-1}\sigma(g)\in K\}\\
&=&\{g\in G(L)/K\mid \sigma(g)=g\}\\
&=&G(F)/G(\mathcal{O}_F),
\end{eqnarray*}
which is the claim.\qed

\subsection{Transitivity of $J\cap P_{\sb}$ on the set of connected components}
For the remainder of Section \ref{secconncomp} we use $M=M_{\sb}$, $P=P_{\sb}$ and $N=N_{\sb}$. 

The strategy of the proof of the second part of Theorem \ref{thmzshk} is as follows: We first show in Proposition \ref{propsurj} that each connected component contains an element of $J\cap \P$ and afterwards connect elements of $J\cap \P\cap \,\ker(\kappa)$ by one-dimensional subvarieties in $X_{\preceq\mu}(b)$.

\begin{prop}\label{propsurj}
Each connected component of $X_{\preceq\mu}(b)$ contains an element of $J\cap \P(L)$.
\end{prop}
The proof of this proposition has two steps: the special case that $b$ is superbasic, and the reduction of the general assertion to this special case. Both the subdivision into these two steps and the proof of the reduction step are inspired by the corresponding proofs for the dimension of affine Deligne-Lusztig varieties, compare \cite{GHKR}, 5. We first consider the case that $b$ is superbasic, i. e. $G=\M=\P$. Then by \cite{GHKR}, 5.9, the root system of $G$ is a disjoint union of root systems of type $A_l$. As we may assume $G$ to be simple and using Lemma \ref{lemnormgm}, we may also assume $G =GL_{h_1}$ for some $h_1> 0$.

Recall the description of $b$ in Section \ref{secnot}. In our context this implies that $b$ is the element of $GL_{h_1}$ that maps $e_i$ to $e_{i+m_1}$ for all $i$ and some $m_1\in\mathbb{Z}$ with $(m_1,h_1)=1$. 

For $\delta\in\mathbb{Z}$ let $\mathbb{A}=\mathbb{A}^1$ if $\delta\neq 0$ and $\mathbb{A}=\mathbb{A}^1\setminus\{-1\}$ if $\delta=0$. Let $\mathcal{GL}_{h_1}$ be the loop group associated to $GL_{h_1}$ over $k$. That is, for every $k$-algebra $R$ we set $\mathcal{GL}_{h_1}(R)=GL_{h_1}(R((t)))$. For $i,\delta\in\mathbb{Z}$ let 
\begin{eqnarray*}
\a_{i,\delta}:\mathbb{A}&\rightarrow &\mathcal{GL}_{h_1}\\
x&\mapsto&\left(e_j\mapsto \begin{cases}
e_j+xe_{j+\delta}&j\equiv i\pmod{h_1}\\
e_j&\text{else}
\end{cases}
\right)
\end{eqnarray*}
Note that $\a_{i,\delta}(x)(te_j)=t\a_{i,\delta}(x)(e_j)$ and that $\a_{i,\delta}=\a_{i+h_1,\delta}$ for all $i$, $\delta$, and $j$.

\begin{lemma}\label{lem1}
Let $g\in GL_{h_1}(L)$. Let $\delta_g\in \mathbb{Z}$ be minimal such that $\a_{i,\delta}(x)\circ g\in gK$ for all $x\in \overline{k}$, $\delta> \delta_g$ and $i$.
\begin{enumerate}
\item $\delta_g\in \{i\in\mathbb{Z}\mid i=-1 \text{ or }i\geq 1\}$.
\item If $\delta_g=-1$, then $gK$ contains an element of $J$.
\item If $\delta_g\geq 0$, then there exists a unique $i_g\in \{1,\dotsc,h_1\}$ with $\a_{i_g,\delta_g}(x)\circ g\notin gK$ for some $x\in\overline{k}$.
\item Let $i,i'\in\{1,\dotsc,h_1\}$, $\delta\geq\delta'$ positive, and $x,x'\in\overline{k}$. Then $[\a_{i,\delta}(x),\a_{i',\delta'}(x')]$ and $\a_{i,\delta}(x)\a_{i,\delta}(x')\a_{i,\delta}(-x-x')$ can be written as products of factors $\a_{i_j,\delta_j}(x_j)$ with $\delta_j>\delta$ (possibly infinite, but converging in the $t$-adic topology). Here, $[.,.]$ denotes the commutator of the two elements.
\end{enumerate}
\end{lemma}

\begin{proof}
We use the bijection between $G(L)/K$ and the set of lattices in $L^{h_1}$ that maps $g\in G(L)$ to $g(\mathcal{O}_L^{h_1})$. Let $\mu=(\mu_j)\in X_*(A)\cong \mathbb{Z}^{h_1}$. A lattice $\Lambda$ corresponds to an element of $X_{\mu}(b)$ if and only if it has a basis $\{x_j\}$ such that the $t^{\mu_j}x_j$ form a basis of $b\sigma(\Lambda)$.

Let $\Lambda\subset L^{h_1}$ be the lattice corresponding to $g\in GL_{h_1}(L)$. As $\det(\a_{i,\delta}(x))\in \mathcal{O}_L^{\times}$ for all $x$, we have $g^{-1}\a_{i,\delta}(x)g\in K$ if and only if it is in $\Mat_{h_1\times h_1}(\mathcal{O}_L)$. This is the case if and only if the following holds: If $v=\sum_{j\in\mathbb{Z}}\beta_je_j\in \Lambda$ with $\beta_j\in\overline{k}$ and $\beta_{i+lh}\neq 0$ for some $l$, then $e_{i+lh+\delta}\in \Lambda$. Indeed, $\a_{i,\delta}(x)\Lambda\subseteq \Lambda$ implies that $\Lambda$ contains both $v$ and $v+\sum_{j\equiv i\pmod{h_1}}\beta_jxe_{j+\delta}$ for some $x\neq 0$.

We assume that $g^{-1}\a_{i,\delta}(x)g\in \Mat_{h_1\times h_1}(\mathcal{O}_L)$ for all $\delta>\delta_0$ and all $i$. Then 
\begin{equation}\label{glvj}v_j=e_j+\sum_{j'>j}\beta_{j'}e_{j'}\in \Lambda
\end{equation} implies $e_l\in \Lambda$ for all $l>j+\delta_0$. Let $j_1$ be minimal such that there is some $v_{j_1}$ as in (\ref{glvj}). Let further $j_2=\max\{j\mid e_j\notin \Lambda\}.$ Then $j_2-j_1=\delta_g$ is minimal with $g^{-1}\a_{i,\delta}(x)g\in \Mat_{h_1\times h_1}(\mathcal{O}_L)$ for all $\delta>\delta_g$ and all $i$. For $\delta=\delta_g$, we have $$g^{-1}\a_{i,\delta}(x)g\in \Mat_{h_1\times h_1}(\mathcal{O}_L)$$ for all $x$ if and only if $i\not\equiv j_1$ modulo $h_1$. 

From this, the first and third assertion follow immediately. For the second note that $\delta_g=j_2-j_1=-1$ implies that $\Lambda=\langle e_i,e_{i+1},\dotsc,e_{i+h_1-1}\rangle_{\mathcal{O}_L}$ for $i=j_1$. Hence $\Lambda=s^i \mathcal{O}_L^{h_1}$ where $s\in J$ with $s(e_j)=e_{j+1}$ for all $j$.

The proof of the last assertion is an easy but tedious computation that is left to the reader.
\end{proof}

\begin{proof}[Proof of Proposition \ref{propsurj} for superbasic $b$]
We may assume that $G=GL_{h_1}$ and that $b$ is of the form discussed above. Let $g\in G(L)$ be a representative of an element of $X_{\preceq\mu}(b)$. We have to show that its connected component contains an element of $J$. For $h_1=1$, each $gK$ contains some $t^i\in J$. From now on we assume that $h_1>1$.

Let $\delta_g$ and $i_g$ be as in Lemma \ref{lem1}. Then by (1) and (2) of this lemma, we may assume that $\delta_g> 0$. Using induction on $\delta_g$ it is enough to show that the connected component of $g$ contains an element $g'$ such that the corresponding number $\delta_{g'}$ is strictly smaller.

We claim that the image of the morphism $\a_{i_g,\delta_g}\circ g:\mathbb{A}\rightarrow X$ is in $X_{\mu}(b)$.

Using Lemma \ref{lem1} (4) and the definition of $\delta_g$ we obtain for every $x\in\overline{k}$
\begin{align*}
Kg^{-1}\a_{i_g,\delta_g}(x)^{-1}b\sigma(\a_{i_g,\delta_g}(x))\sigma(g)K&=Kg^{-1}\a_{i_g,\delta_g}(-x)\a_{i_g+m_1,\delta_g}(\sigma(x))b\sigma(g)K\\
&=Kg^{-1}\a_{i_g+m_1,\delta_g}(\sigma(x))\a_{i_g,\delta_g}(-x)b\sigma(g)K\\
&=K g^{-1}\a_{i_g+m_1,\delta_g}(\sigma(x))b\a_{i_g-m_1,\delta_g}(-x)\sigma(g)K.\\
\intertext{As $h_1>1$, neither $i_g+m_1$ nor $i_g-m_1$ is congruent to $i_g$ modulo $h_1$. Together with the uniqueness of $i_g$ we obtain}
&= Kg^{-1}b\sigma(g)K\\
&=Kt^{\mu}K.
\end{align*}
Hence $\a_{i_g,\delta_g}(x)\circ g\in X_{\mu}(b)$.

As $X_{\preceq \mu}(b)$ satisfies the valuation criterion of properness, this morphism induces a morphism $\phi:\mathbb{P}^1\rightarrow X_{\preceq\mu}(b)$. For $x=[1:0]$, we have $\phi(x)=g$. Let $g'=\phi([0:1])$. Then $g$ and $g'$ are in the same connected component of $X_{\preceq\mu}(b)$. It remains to show that $\delta_{g'}<\delta_g$. For every $x\in\mathbb{A}$ we have to consider $\a_{i,\delta}(x)\circ g'$ for all $\delta\geq\delta_g>0$. It is the image of $[0:1]$ in the closure of the family $\a_{i,\delta}(x)\circ \a_{i_g,\delta_g}(y)\circ g\in X$ with $y\in \mathbb{A}\subseteq \mathbb{P}^1$. Lemma \ref{lem1} (4) shows that
\begin{equation}\label{glhurra}
\a_{i,\delta}(x)\circ \a_{i_g,\delta_g}(y)\circ gK=\a_{i_g,\delta_g}(y)\circ \a_{i,\delta}(x)\circ gK
\end{equation} 
for all $\delta\geq\delta_g$. But for $\delta>\delta_g$ or $\delta=\delta_g$ and $i\neq i_g$ this is equal to $\a_{i_g,\delta_g}(y)\circ gK$. Thus the image of $[0:1]$ is again $g'$, and we obtain $\delta_{g'}\leq \delta_g$. Consider (\ref{glhurra}) for $\delta=\delta_g$ and $i=i_g$. By Lemma \ref{lem1} (4) the right hand side is equal to $$\a_{i_g,\delta_g}(x+y)\circ gK.$$ The image of $[0:1]$ in this family is again $g'$. Hence $\delta_{g'}<\delta_g$, which completes the induction step.
\end{proof}

\begin{proof}[Proof of Proposition \ref{propsurj}]
It remains to show that for every $g\in X_{\mu}(b)$ there is a $g'$ in the same connected component of $X_{\preceq\mu}(b)$ with $g'=jm$ where $j\in J\cap \P(L)$ and $m\in \M(L)$. In the following we find such data which even satisfy $j\in J\cap N(L)$.

By the Iwasawa decomposition we write $g=nmk$ with $n\in N(L)$, $m\in M(L)$, and $k\in K$. We may thus assume that $g\in gK$ is already equal to $nm$. Then 
\begin{equation*}
g^{-1}b\sigma(g)=m^{-1}\left(n^{-1}b\sigma(n)b^{-1}\right)b\sigma(m).
\end{equation*}
We abbreviate the expression in the bracket, which is in $N(L)$, by $\tilde{n}$. Let $\mathcal{N}$ be the loop group associated to $N$ over $\overline{k}$ and
\begin{eqnarray*}
f_b:\mathcal{N}&\rightarrow&\mathcal{N}\\
n_0&\mapsto&n_0^{-1}b\sigma(n_0)b^{-1}.
\end{eqnarray*}
This is a morphism of ind-schemes over $\Spec(\overline{k})$. We have $\tilde{n}=f_b(n)$.

Let $\chi\in X_*(A_P)$ be such that $\langle\alpha,\chi\rangle>0$ for every simple root $\alpha$ of $A$ in $N$. Let
\begin{eqnarray*}
\phi: \mathbb{A}_{\overline{k}}^1\setminus\{0\} &\rightarrow &\mathcal{N}\\
x&\mapsto& x^{\chi}\tilde{n}x^{-\chi}.
\end{eqnarray*}
Let $\alpha$ be a root of $A$ in $N$ and let $U_{\alpha}$ denote the corresponding root subgroup. Conjugation by $x^{\chi}$ maps $U_{\alpha}(y)$ to $U_{\alpha}(x^jy)$ where $j=\langle\alpha,\chi\rangle>0$. Especially, $\phi$ has an extension to a morphism $\phi:\mathbb{A}^1\rightarrow \mathcal{N}$ that maps $0$ to $\id$. Besides, $$m^{-1}\phi(x)b\sigma(m)=x^{\chi}m^{-1}\tilde{n}b\sigma(m)x^{-\chi}\in Kt^{\mu}K$$ for every $x\neq 0$. This implies that also $m^{-1}\phi(0)b\sigma(m)\in Kt^{\mu'}K$ for some $\mu'\preceq\mu$. 

We construct in the following a finite \'{e}tale morphism $p:\Spec(R)\rightarrow\mathbb{A}^1$ for some $R$ such that the restriction of $p$ to each connected component of $\Spec(R)$ is surjective. Besides we construct a morphism $\psi:\Spec(R)\rightarrow \mathcal{N}$ such that $(f_b\circ \psi(x))b\sigma(m)K= (\phi\circ p(x))b\sigma(m)K$ for all geometric points $x$ of $\Spec(R)$. We also require that there is a geometric point $x_1$ with $\psi(x_1)m=g$ as elements of $X$. The first condition implies that for every $x$ there is some $\mu'\preceq\mu$ such that $m^{-1}f_b(\psi(x))b\sigma(m)\in Kt^{\mu'}K$. Especially, $\psi(x)m\in X_{\preceq\mu}(b)$. Let $x_0$ be in the same connected component of $\Spec(R)$ as $x_1$ and with $p(x_0)=0$. Then $\psi(x_0)\in J$ and $\psi(x_0)m$ is in the same connected component of $X_{\preceq\mu}(b)$ as $g$. 

Let $\alpha_i$ with $i=1,\dotsc,i_0$ be the roots of $A$ in $N$. Write each $\alpha_i$ as a linear combination of simple roots of $G$. Let $N[j]\subseteq N$ for $j\geq 1$ be the normal subgroup that is generated by the $\alpha_i$ such that the sum of the corresponding coefficients of the simple roots of $A$ in $N$ is at least $j$. Then the commutator of elements of $N[j]$ and $N[j']$ is in $N[j+j']$. Especially, $N[j]/N[j+1]$ is abelian. Let $N=N_1\supset N_2\supset \dotsm$ be a refinement of $N=N[1]\supset N[2]\supset \dotsm$ with $N_i/N_{i+1}\cong \mathbb{G}_a$. We choose $N_i$ and the ordering of the $\alpha_i$ such that $N_i$ is generated by all $U_{\alpha_{i'}}$ with $i'\geq i$. We will use induction on $i$ to construct $R_i$, $p_i$ and $\psi_i$ that satisfy the first of the claimed conditions up to factors in $N_{i+1}$. More precisely, we want $(\phi\circ p_i(x))b\sigma(m)\in (f_b\circ\psi(x))N_{i+1}(R_{i}((t)))b\sigma(m)K$. As only finitely many of the $N_{i}$ are non-trivial, we then obtain $R, p$ and $\psi$ with $(f_b\circ \psi(x))b\sigma(m)K= (\phi\circ p(x))b\sigma(m)K$. The existence of $x_1$ is shown afterwards.

{\bf Claim.} {\it Let $R_{i-1}$ be an \'{e}tale extension of $\overline{k}[s]$ such that each connected component of $\Spec(R_{i-1})$ maps surjectively to $\Spec(\overline{k}[s])$ and let $y_i\in N_i(R_{i-1}((t)))$. Then we can write $y_i\in U_{\alpha_i}(\beta_i) N_{i+1}(R_{i-1}((t)))$ for some $\beta_i\in R_{i-1}((t))$. Let $n(i)\in\mathbb{N}$. Let $j$ be minimal with $N_{i}\subseteq N[j]$. Then there exists a finite \'{e}tale extension $R_{i}$ of $R_{i-1}$ with the same surjectivity property and $z\in N[j](R_{i}((t)))$ such that $f_b(z)= U_{\alpha_i}(\beta_i+\varepsilon_i)$ with $\varepsilon_{i}\in   t^{n(i)}R_{i}[[t]]$.}

This claim implies the induction step. Indeed, let $y_{\phi}\in N(\overline{k}[s]((t)))$ be the element associated to $\phi$. Assume that we constructed $z_{i-1}\in N(R_{i-1}((t)))$ with $$y_{\phi}b\sigma(m)\in f_b(z_{i-1})\delta_ib\sigma(m)K$$ for some $\delta_i\in N_{i}(R_{i-1}((t)))$. Set $z_{i}=z_{i-1}\tilde{z}\in N(R_{i}((t)))$ where we have to define $R_{i}$ and $\tilde{z}\in N[j](R_{i}((t)))$ such that $y_{\phi}b\sigma(m)\in f_b(z_{i})N_{i+1}(R_{i}((t)))b\sigma(m)K$. We have
\begin{equation*}
f_b(z_{i})=\tilde{z}^{-1}f_{b}(z_{i-1})(b\sigma(\tilde{z})b^{-1}).
\end{equation*}
As $\tilde{z}, b\sigma(\tilde{z})b^{-1}\in N[j](R_{i}((t)))$, their commutators with $f_b(z_{i-1})$ are in $N[j+1](R_{i}((t)))\subseteq N_{i+1}(R_{i}((t)))$. Hence, modulo $N_{i+1}$,
\begin{eqnarray*}
f_b(z_{i})&\equiv&f_b(z_{i-1})f_b(\tilde{z}).
\end{eqnarray*}
Using the claim for sufficiently large $n(i)$ and $y_i=\delta_i$, one finds $\tilde{z}$ such that $$f_b(\tilde{z})^{-1}\delta_{i}b\sigma(m)\in N_{i+1}(R_{i}((t)))b\sigma(m)K,$$ which is what we wanted.
 
To prove the claim, recall the description of $b$ as $t^{\mu_{M}}w$ with $w\in W_M$ from Section \ref{secnot}. There is some $l_0>0$ such that $w^{l_0}=\id$. Then $$(b\sigma)^{l_0}U_{\alpha_i}(\beta_i)b^{-l_0}=U_{\alpha_i}(t^{l_0\langle\alpha_i,\nu\rangle}\sigma^{l_0}(\beta_i)).$$

We distinguish two cases.

{\bf Case 1.} $\langle\alpha_i,\nu\rangle>0$

Here, let $R_{i+1}=R_i$. As $\langle\alpha_i,\nu\rangle>0$, the $(b\sigma)^lU_{\alpha_i}(\beta_i)b^{-l}$ converge to $\id$ for large $l$. For $l\in\mathbb{N}$ let $$z(l)=U_{\alpha_i}(-\beta_i)(bU_{\alpha_i}(\sigma(-\beta_i))b^{-1})\dotsm(b^lU_{\alpha_i}(\sigma^l(-\beta_i))b^{-l})$$ and let $z$ be the limit of this sequence. Then 
\begin{eqnarray*}
f_b(z(l))&=&U_{\alpha_i}(\beta_i)\dotsm(b^lU_{\alpha_i}(\sigma^l(\beta_i))b^{-l})(b^{l+1}U_{\alpha_i}(\sigma^{l+1}(-\beta_i))b^{-l-1})\dotsm (bU_{\alpha_i}(\sigma(-\beta_i))b^{-1})\\
&=&U_{\alpha_i}(\beta_i)[(bU_{\alpha_i}(\sigma(\beta_i))b^{-1})\dotsm(b^lU_{\alpha_i}(\sigma^l(\beta_i))b^{-l}),b^{l+1}U_{\alpha_i}(\sigma^{l+1}(-\beta_i))b^{-l-1}]\\&&\times(b^{l+1}U_{\alpha_i}(\sigma^{l+1}(-\beta_i))b^{-l-1}).
\end{eqnarray*}
Here, $[.,.]$ denotes the commutator of the two elements. Thus, $f_b(z)=U_{\alpha_i}(\beta_i)$.

{\bf Case 2.} $\langle\alpha_i,\nu\rangle=0$

In this case, the $b^lU_{\alpha_i}(\beta_i)b^{-l}$ do not converge to $\id$. There is some $l_0>0$ such that $$(b\sigma)^{l_0}U_{\alpha_i}(\beta_i)b^{-l_0}=U_{\alpha_i}(\sigma^{l_0}(\beta_i)).$$ Let $l_0>0$ be minimal with this property. We set 
\begin{equation}\label{glz}
z=U_{\alpha_i}(z_0)(b\sigma(U_{\alpha_i}(z_0))b^{-1})\dotsm ((b\sigma)^{l_0-1}U_{\alpha_i}(z_0)b^{-l_0+1})
\end{equation}
with $z_0\in R_{i+1}((t))$ where $R_{i+1}$ is yet to define. Let $j$ be minimal with $U_{\alpha_i}\subseteq N[j]$. As $b\in M(L)$, it is also minimal with $b^lU_{\alpha_i}b^{-l}\subseteq N[j]$ for each $l$. Then up to factors $U_{\alpha_{i'}}(z')$ in $N[j+1]\subseteq N_{i+1}$ we obtain
\begin{eqnarray*}
z^{-1}b\sigma(z)b^{-1}&=&z^{-1}U_{\alpha_i}(-z_0)z (b\sigma)^{l_0}U_{\alpha_i}(z_0)b^{-l_0}\\
&\equiv&U_{\alpha_i}(\sigma^{l_0}(z_0)-z_0).
\end{eqnarray*}
Thus $z_0$ has to satisfy the equation
$$\sigma^{l_0}(z_0)\equiv z_0+\beta_i \pmod{t^{n(i)}R_{i+1}[[t]]}.$$
We write $\beta_i=\sum_{j\geq j_0} \beta_{i,j}t^j$ with $\beta_{i,j}\in R_i$. We choose some $a_i\leq j_0$, and let $\beta_{i,j}=0$ for $j<j_0$. Let $z_0=\sum_{j= a}^{n(i)-1} z_{0,j}t^j$. Let $R_{i+1}$ be the \'{e}tale extension of $R_i$ defined by adjoining a root of the polynomial $z_{0,j}^{q^{l_0}}-z_{0,j}-\beta_{i,j}$ for every $j$ with $a_i\leq j<n(i)$. Then $z$ as in (\ref{glz}) satisfies all required properties. The surjectivity property follows from the explicit description of the extension $R_{i+1}$ of $R_i$.

It remains to show the existence of a geometric point $x_1$ of $\Spec(R)$ over $1\in\mathbb{A}^1$ with $\psi(x_1)m=g$ as elements of $X$. Note that for any $a\in\mathbb{Z}$, there are only finitely many cosets $n'mK$ of elements $n'\in N(t^a\mathcal{O}_L)$ with $f_b(n')=f_b(n)$. Here $N(t^a\mathcal{O}_L)$ denotes the subgroup of $N$ generated by all $U_{\alpha}(t^a\mathcal{O}_L)$ for $\alpha$ positive. Indeed, the elements $n'm$ with $f_b(n')=f_b(n)$ are in the $J$-orbit of $g$. Taking $K$-cosets together with the restriction to be in $N(t^a\mathcal{O}_L)$, we obtain a finite number of them. The construction of $R$ shows that for some fixed $a$ we can choose the $a_i$ in the second step of the claim small enough such that the following holds: For each of the cosets $n'mK$ above there is a representative $n'=\psi(x')$ for some $x'\in\Spec(R)$ lying over $1\in\mathbb{A}^1$. Especially, for sufficiently small $a_i$, a point $x_1$ as above exists. 
\end{proof}

The following corollary is a generalization of \cite{KottRapo}, Proposition 4.8. There, a similar statement is shown for $G=GL_n$ and $G=GSp_{2n}$.
\begin{kor}\label{corsurj}
Let $b\in G(L)$ and $\nu\in X_*(A)$ its Newton point. Let $\mu_{\min}$ be as in Section \ref{secnot}. Then $\mu_{\min}$ is minimal among the dominant elements $\mu$ with $X_{\mu}(b)\neq\emptyset$. More precisely, for every dominant $\mu\in X_*(A)$ with $\nu\preceq\mu$ and $\kappa_G(b)=\mu$ (as elements of $\pi_1(G)$) we have $\nu\preceq\mu_{\min}\preceq\mu$.
\end{kor}
\begin{proof}
The assumptions on $\mu$ are equivalent to $X_{\preceq\mu}(b)\neq \emptyset$. Then the proof of Proposition \ref{propsurj} shows that $X_{\preceq\mu}(b)$ contains an element of $J$. Our choice of the normal form of $b$ in Section \ref{secnot} implies that the image of $J$ in $X$ lies in $X_{\mu_{\min}}(b)$. Thus $\mu_{\min}\preceq\mu$.
\end{proof}

\subsection{Proof of Theorem \ref{thmzshk} (2)}
To deduce Theorem \ref{thmzshk} (2) from Proposition \ref{propsurj}, we also need the following proposition.

\begin{prop}\label{props0}
Let $\alpha_0$ be a simple root of $A$ in $N$ and let the assumptions of Theorem \ref{thmzshk} (2) be satisfied. Then there exists a root $\alpha$ with $\alpha_0^{\vee}=\alpha^{\vee}$ as elements of $\pi_1(M)$ and such that $t^{-\alpha^{\vee}}$ lies in the connected component of the identity of $X_{\preceq\mu}(b)$. 
\end{prop}

\begin{proof}[Proof of Theorem \ref{thmzshk} (2)]
By Proposition \ref{propsurj}, the $J\cap \P$-equivariant map $\varphi:(J\cap \P)/K\rightarrow \pi_0(X_{\preceq\mu}(b))$ induced by the inclusion of the left hand side into $X_{\preceq\mu}(b)$ is surjective. We have to show that $\ker(\kappa_G)\cap J\cap \P(L)$ is mapped to the connected component of the identity. 

Let $j\in J\cap \P$ and $j=mn$ with $m\in \M(L)$ and $n\in \N(L)$. Then 
\begin{equation*}
b=j^{-1}b\sigma(j)= \left(n^{-1}m^{-1}b\sigma(mn) (m^{-1}b\sigma(m))^{-1}\right)\left(m^{-1}b\sigma(m)\right)
\end{equation*}
As the first bracket is in $\N$ and the second in $\M$, we obtain that $m^{-1}b\sigma(m)=b$. Hence $m\in J$ which implies that also $n\in J$. Note that by the definition of $\kappa$, unipotent elements are in the kernel of $\kappa$. Hence $j\in\ker(\kappa)$ if and only if $m$ and $n$ are in $\ker(\kappa)$. As $\varphi$ is $J$-equivariant, it is thus enough to consider elements of $J\cap \M$ and $J\cap \N$ separately.

Let $j\in J\cap \M\cap \ker(\kappa_M)$. From the proof of Proposition \ref{propsurj} for superbasic $b$ we obtain that $j$ is in the same connected component as $\id$. Thus the restriction of $\varphi$ to $J\cap M$ factors through a quotient of $\pi_1(M)$. Besides, Proposition \ref{props0} implies that for every simple root $\alpha_0$ of $A$ in $N$, there is some $g\in M(L)$ in the connected component of $\id$ in $X_{\preceq\mu}(b)$ with $\kappa_M(g)=-\alpha_0^{\vee}$. Hence, the restriction of $\varphi$ to $J\cap M$ even factors through $\pi_1(G)$.

It remains to consider elements $j\in J\cap \N$. Each such $j$ can be written as a finite product of $U_{\alpha_i}(x_i)$ with $\alpha_i$ a root of $A$ in $N$ and $x_i\in L$. There is a $y\in X_*(A_P)\subseteq X_*(A)$ with $y$ trivial in $\pi_1(G)$ and such that $\langle \alpha_i,y\rangle\geq -v_t(x_i)$ for all $i$. Then $j$ is conjugate via $t^y\in A(L)\cap J$ to an element $\tilde{j}$ of $\N(\mathcal{O}_L)\cap J\subseteq K$. We have seen above that $\varphi$ maps $t^y\in J\cap M(F)$ to the connected component of the identity. As $\tilde{j}\in J\cap K$, it is also mapped to this component. Hence, $j$ is in the connected component of the identity of $X_{\preceq \mu}(b)$.
\end{proof}

\begin{proof}[Proof of Proposition \ref{props0}]
We first show that it is enough to prove that there is an $\alpha$ with $\alpha^{\vee}=\alpha_0^{\vee}$ in $\pi_1(M)$ such that the morphism $g:\mathbb{A}^1\rightarrow X$ which maps $x$ to $U_{\alpha}(t^{-1}x)$ factors through $X_{\preceq\mu}(b)$. If this is the case, fact that $X_{\preceq\mu}(b)$ satisfies the valuation criterion of properness leads to an extension $g:\mathbb{P}^1\rightarrow X_{\preceq\mu}(b)$. Thus $g([1:0])=\id$ and $g([0:1])$ are in the same connected component. But for $x\neq [1:0]$, an $SL_2$-calculation shows that 
\begin{equation}\label{glU+-}
U_{\alpha}(t^{-1}x)= U_{-\alpha}(tx^{-1})t^{-\alpha^{\vee}}k
\end{equation}
for some $k\in K$. Thus $g([0:1])=t^{-\alpha^{\vee}}$ as elements of $X$.

By passing to a suitable subgroup of $G$, we may assume that the Dynkin diagram consisting of the roots of $M$, and of $\alpha_0$, is connected. 

Recall the description of $b$ as $t^{\mu_M}w$ with $w=(w_i)\in W_M$ and $\mu_M\in X_*(A)$ from Section \ref{secnot}. We have to show that for all $x$ there exists some $\mu'\preceq\mu$ such that
\begin{equation}\label{glmain}
A:=g(x)^{-1}b\sigma(g(x))=U_{\alpha}(-t^{-1}x)\left( b U_{\alpha}(t^{-1}\sigma(x))b^{-1}\right) b\in Kt^{\mu'}K.
\end{equation}
We have
\begin{eqnarray}
\nonumber bg(\sigma(x))b^{-1}&=&t^{\mu_M}U_{w(\alpha)}(t^{-1}\sigma(x))t^{-\mu_M}\\
\label{glmain2}&=&U_{w(\alpha)}(t^{\langle w(\alpha),\mu_M\rangle-1}\sigma(x))
\end{eqnarray}
where $w(\alpha)$ is such that $wU_{\alpha}w^{-1}=U_{w(\alpha)}$.

We consider several cases. Essentially, we have to distinguish between the different root systems and values of $\langle\alpha_0,\mu_M\rangle$. However, many of these cases are treated similarly, especially, if the length of $\alpha_0$ is the same.

{\bf Case 1. $\alpha_0$ is the only root of $G$.} 

In this case we may assume $G=GL_2$. Besides, $w=\id$. Let $\alpha=\alpha_0$. Recall that $\mu_M=\nu$ is dominant, so $\langle\alpha,\mu_M\rangle\geq 0$. Equations (\ref{glmain}) and (\ref{glmain2}) show that $A=U_{\alpha}(t^{-1}y)t^{\mu_M}$ for some $y\in \mathcal{O}_L$. If $y\in t\mathcal{O}_L$ the claim follows. Assume $y\in \mathcal{O}_L^{\times}$. Together with (\ref{glU+-}) we obtain that 
\begin{align*}
U_{\alpha}(t^{-1}y)t^{\mu_M}&\in K t^{\alpha^{\vee}}U_{-\alpha}(ty^{-1})t^{\mu_M}K\\
&=Kt^{\mu_M+\alpha^{\vee}}U_{-\alpha}(t^{1+\langle\alpha,\mu_M\rangle}y^{-1})K\\
&= Kt^{\mu_M+\alpha^{\vee}}K.
\end{align*}
As $\alpha$ is the only simple root, $\mu_M+\alpha^{\vee}$ is the minimal dominant $\mu$ with $\mu_M\preceq\mu$ and $\mu_M\neq\mu$. 

{\bf Case 2. All roots of $G$ have the same length.} 

As the root system of $M$ is a sum of root systems of type $A_l$, this condition is equivalent to the condition that the neighbours of $\alpha_0$ in the Dynkin diagram are of the same length as $\alpha_0$.

Again we choose $\alpha=\alpha_0$. The first step is to compute $w(\alpha)$. We describe the effect of each of the parts $w_i$ of $w$ separately. Let $\alpha^{(i)}$ be the image of $\alpha$ in $X_*(A_i)\cong \mathbb{Z}^{h_i}/\mathbb{Z}$. Let $e_j$ for $1\leq j \leq h_i$ be the standard system of generators of $X_*(A_i)$. We may assume that the simple roots of $M_i$ are the $\beta_j=e_j-e_{j+1}$. Let $\beta_{j_0}$ be the neighbour of $\alpha$ in the Dynkin diagram of $G$. Note that the classification of Dynkin diagrams shows that $j_0\in\{1,2,3,h_i-3,h_i-2,h_i-1\}$. For symmetry reasons it is enough to consider the first three cases. Thus $\alpha^{(i)}$ is equal to $-e_1$, $-e_1-e_2$, or $-e_1-e_2-e_3$. Hence $w_i(\alpha^{(i)})$ is equal to $-e_{1+m_i}$, $-e_{1+m_i}-e_{2+m_i}$, or $-e_{1+m_i}-e_{2+m_i}-e_{3+m_i}$. Here we identify $e_{j}$ with $e_{j-h_i}$ if $j>h_i$. Recall that $\mu_{M_i}=(1,\dotsc,1,0\dotsc,0)$ with multiplicities $m_i$ and $h_i-m_i$. Thus
\begin{equation*}
\langle w_i(\alpha)-\alpha,\mu_M\rangle=\langle \alpha-w_i(\alpha),\alpha^{\vee}\rangle
\end{equation*}
for all $i$. This implies
\begin{eqnarray}\label{glgain2}
\langle w(\alpha)-\alpha,\mu_M+\alpha^{\vee}\rangle&=&\sum_i\langle w_i\dotsm w_1(\alpha)-w_{i-1}\dotsm w_1(\alpha),\mu_M+\alpha^{\vee}\rangle\\
\nonumber &=&\sum_i\langle w_i(\alpha)-\alpha,\mu_M+\alpha^{\vee}\rangle\\
\nonumber &=&0
\end{eqnarray}
This explicit description also shows that $w(\alpha)$ is a positive root. Besides, with $$\nu_{M_i}=(m_i/h_i,\dotsc,m_i/h_i)$$ we obtain that 
\begin{equation}\label{eqlast}
\langle w(\alpha),\mu_M\rangle>\langle w(\alpha),\nu\rangle\geq 0.
\end{equation}
We consider the case $\langle w(\alpha),\alpha^{\vee}\rangle\geq 0$. Then $\langle a\alpha+bw(\alpha),\alpha^{\vee}\rangle\geq 2a$ and $\langle a\alpha+bw(\alpha),w(\alpha)^{\vee}\rangle\geq 2b$. Hence for all $a,b>0$, the linear combination $a\alpha+bw(\alpha)$ cannot be a root. This implies that $U_{\alpha}$ and $U_{w(\alpha)}$ commute. From (\ref{glmain}) and (\ref{glmain2}), we obtain for $\langle\alpha,\mu_M\rangle\geq 0$
\begin{align*}
A&\in KU_{\alpha}(-t^{-1}x)U_{w(\alpha)}(y)t^{\mu_M}K\\
&=KU_{w(\alpha)}(y)U_{\alpha}(-t^{-1}x)t^{\mu_M}K.\\
\intertext{for some $y\in\mathcal{O}_L$. By (\ref{glU+-}) this is}
&=Kt^{\alpha^{\vee}}U_{-\alpha}(-tx^{-1})t^{\mu_M}K\\
&=Kt^{\alpha^{\vee}+\mu_M}K.
\end{align*}
It only remains to show $(\mu_M+\alpha^{\vee})_{\dom}\preceq\mu$.

For $\langle\alpha,\mu_M\rangle< 0$ we obtain 
\begin{eqnarray*}
A& \in &KU_{\alpha}(-t^{-1}x)U_{w(\alpha)}(y))t^{\mu_M}K\\
&=&KU_{\alpha}(-t^{-1}x)t^{\mu_M}K\\
&=&Kt^{\mu_M}K
\end{eqnarray*}
for some $y\in\mathcal{O}_L$. Here, the last equation follows from (\ref{eqlast}).

Let now $\langle w(\alpha),\alpha^{\vee}\rangle=-1$. Then we similarly obtain that $\alpha+w(\alpha)$ is the only positive linear combination that can be a root. Besides, we know from (\ref{glgain2}) that $\langle w(\alpha)-\alpha,\mu_M\rangle= 3$. If in this case $\langle w(\alpha),\mu_M\rangle\geq 2$, (\ref{glmain}) and (\ref{glmain2}) imply
\begin{align*}
A&\in U_{\alpha}(-xt^{-1})U_{w(\alpha)}(ty)t^{\mu_M}K\\
&\subseteq  KU_{\alpha+w(\alpha)}(-cxy)U_{w(\alpha)}(ty)U_{\alpha}(-xt^{-1})t^{\mu_M}K\\
&= KU_{\alpha}(-xt^{-1})t^{\mu_M}K\\
\intertext{for some $y\in\mathcal{O}_L$ and $c\in k$. Here we set $U_{\alpha+w(\alpha)}=1$ if $\alpha+w(\alpha)$ is not a root. For $\langle\alpha,\mu_M\rangle\geq 0$ we use (\ref{glU+-}) and obtain}
&= \begin{cases}
Kt^{\alpha^{\vee}+\mu_M}K&\text{if }\langle\alpha,\mu_M\rangle\geq 0\\
Kt^{\mu_M}K&\text{else. }
\end{cases}
\end{align*}
The calculation for the remaining case $\langle \alpha,\mu_M\rangle\leq -2$ is similar.

If $\langle\alpha,\mu_M\rangle\geq 0$, it remains to compare $(\mu_M+\alpha^{\vee})_{\dom}$ and the minimal dominant $\mu$ with $\mu_{\min}\preceq \mu$ and $\mu_{\min}\neq\mu$ in $\pi_1(M)$. We describe how to compute $(\mu_M+\alpha^{\vee})_{\dom}$ and $\mu$. For both Hodge vectors one can start with $\mu_M+\alpha^{\vee}$ and then successively add coroots until one reaches the desired element. For $(\mu_M+\alpha^{\vee})_{\dom}$ one replaces $\mu_M+\alpha^{\vee}$ by images $\mu_i$ under reflections in the Weyl group. In the $i$th step one adds $-\beta_i^{\vee}\frac{\langle\beta_i,\mu_{i-1}\rangle}{\langle\beta_i,\beta_i^{\vee}\rangle}$ to $\mu_{i-1}$ for some simple root $\beta_i$. One chooses the roots $\beta_i$ such that $\langle\beta_i,\mu_{i-1}\rangle<0$. Then one can inductively show that in each step $-\langle\beta_i,\mu_{i-1}\rangle/\langle\beta_i,\beta_i^{\vee}\rangle= 1$. But this is also exactly the description of how to obtain $\mu$ from $\mu_{\min}+\alpha^{\vee}$. Thus the two Hodge vectors are equal.

{\bf Case 3. $\alpha_0$ is a long root of $G$.}
In this case let $\alpha=\alpha_0$. As $w(\alpha)$ is in the $W$-orbit of $\alpha$, it is also a long root. The same arguments as in the preceding case show that (\ref{glgain2}) holds in this case, that $\langle w(\alpha),\mu_M\rangle>0$, and that $(\mu_M+\alpha^{\vee})_{\dom}\preceq\mu$. If $\alpha_i$ is the shorter neighbour of $\alpha$, then $\langle w_i(\alpha)-\alpha,\mu_M\rangle\geq 2$. The same calculations and subcases as in Case 2 show the claim.

{\bf Case 4. $\alpha_0$ is a short root of $G$.} 

Let $\beta$ be the longer neighbour of $\alpha_0$ in the Dynkin diagram and $\alpha=(\alpha_0^{\vee}+\beta^{\vee})^{\vee}$. Then $\alpha$ is a long root of $G$. As in Case 2 one can show that $(\mu_{\min}+\alpha^{\vee})_{\dom}\preceq\mu$. Note that the corresponding assertion for $\alpha_0$ instead of $\alpha$ does not hold in general in this case. We treat the different possible root systems separately.

If $G$ is of type $G_2$, the long roots $\beta$ and $\alpha=3\alpha_0+\beta$ generate a sub-root system of type $A_2$, and the corresponding subgroup of $G$ contains $b$. Thus the assertion for this case follows from Case 2.

If $G$ is of type $C_n$ with $n\geq 2$, the roots of $M$ and $\alpha=2\alpha_0+\beta$ generate a sub-root system of type $C_{n-1}+A_1$. The corresponding Levi subgroup of $G$ contains $b$, and $\alpha$ is a long root. This case is covered by Case 3.

If $G$ is of type $B_n$ with $n\geq 3$, again $\alpha=2\alpha_0+\beta$. We have $n=h_1$ where $h_1$ is as in the definition of the normal form of $b$ in Section \ref{secnot}. An explicit description of the root system is given by the simple roots $\beta=e_1-e_2,\dotsc,e_{n-1}-e_n$, and $\alpha_0=-e_1$. Then $\alpha=-e_1-e_2$ and $w(\alpha)=-e_{m_1+1}-e_{m_1+2}$. The sum of the coefficients of each root with respect to the basis $\{e_i\}$ is at least $-2$. Hence no linear combination $a\alpha+b w(\alpha)$ with $a,b>0$ is a root. Again, the image of $\mu$ in $X_*(A_1)$ is equal to $(c+1,\dotsc,c+1,c,\dotsc,c)$ with multiplicities $m_1$ and $h_1-m_1$ and $c\in \mathbb{Z}$. As $h_1=n>2$, we have $\langle w(\alpha),\mu_M\rangle >\langle w(\alpha),\nu\rangle\geq 0$. Thus we obtain
\begin{eqnarray*}
A&=&U_{\alpha}(-xt^{-1})U_{w(\alpha)}(y)t^{\mu_M}w\\
&\in& KU_{\alpha}(-xt^{-1})t^{\mu_M}K
\end{eqnarray*}
for some $y\in \mathcal{O}_L$. If $\langle\alpha,\mu_M\rangle <0$, this is in $Kt^{\mu_M}K$. Else we use again (\ref{glU+-}) and obtain $U_{\alpha}(-xt^{-1})t^{\mu_M}\in Kt^{\mu_M+\alpha^{\vee}}K$.

Let now $G$ be of type $F_4$. Again $\alpha=2\alpha_0+\beta$. Using the explicit description of this root system one sees that $\alpha+w(\alpha)$ is the only positive linear combination of $\alpha$ and $w(\alpha)$ that can be a root. Besides, the explicit calculation of $w(\alpha)$ yields $\langle w(\alpha)-\alpha,\mu_M\rangle\geq 3$. Using the same calculation as in Case 2 for $\langle w(\alpha),\alpha^{\vee}\rangle=-1$, the claim follows. 
\end{proof}

\section{Non-closed affine Deligne-Lusztig varieties}\label{secopen}
We give an example of a non-closed affine Deligne-Lusztig variety for $G=GL_5$ with $\pi_0(X_{\mu}(b))\neq \pi_0(X_{\preceq\mu}(b))$. We use the interpretation of $X$ as a set of lattices in $L^5$ that is introduced in the proof of Lemma \ref{lem1}.

Let $B\subseteq GL_5$ be the Borel subgroup of upper triangular matrices and let $A$ be the diagonal torus.

We denote a basis of $L^{5}$ by $e_{1,0}$, $e_{1,1}$, $e_{2,0}$, $e_{2,1}$, and $e_{2,2}$. For $i\in \mathbb{Z}$ let $e_{1,i}=te_{1,i-2}$ and $e_{2,i}=te_{2,i-3}$. Let $b\in G(L)$ with $b\sigma(e_{i,j})=e_{i,j+1}$. Then $\nu=(\frac{1}{2},\frac{1}{2},\frac{1}{3},\frac{1}{3},\frac{1}{3})\in\mathbb{Q}^5\cong X_*(A)_{\mathbb{Q}}$. Let further $\mu=(2,0,0,0,0)$. Then Theorem \ref{thmzshk} implies $\pi_0(X_{\preceq\mu}(b))\cong \mathbb{Z}$. We prove in this section that 
\begin{equation}
\pi_0(X_{\mu}(b))\cong \mathbb{Z}^2 \cong J/(K\cap J).
\end{equation}
More precisely, we define a morphism $\coprod_{\mathbb{Z}^2}\mathbb{A}^2\rightarrow X_{\mu}(b)$ which is a bijection on $\overline{k}$-valued points.

Let $\Lambda\subset L^5$ be the lattice corresponding to a point of $X_{\mu}(b)$. As there is only one $\mu_i\neq 0$, we have $\Lambda/(b\sigma(\Lambda)+t\Lambda)\cong \overline{k}$. Each $v\in \Lambda\setminus (b\sigma(\Lambda)+t\Lambda)$ generates $\Lambda$ as a $b\sigma$-invariant $\mathcal{O}_L$-module. We renormalize the second indices of the $e_{i,j}$ by a suitable shift such that $v=\sum_{i\in\{1,2\},j\geq 0} \beta_{i,j}e_{i,j}$ with $\beta_{i,j}\in\overline{k}$ and $\beta_{i,0}\neq 0$. By multiplying with $\beta_{1,0}^{-1}\in\overline{k}^{\times}$ we may assume that $\beta_{1,0}=1$. We have $(b\sigma)^lv=e_{1,l}+\sum_{j> 0} \beta_{1,j}^{\sigma^l}e_{1,j+l}+\sum_{j\geq 0} \beta_{2,j}^{\sigma^l}e_{2,j+l}$, hence we can modify $v$ by an infinite but converging linear combination of these elements for $l>0$ to obtain an element of $\Lambda\setminus (b\sigma(\Lambda)+t\Lambda)$ of the form $e_{1,0}+\sum_{j\geq 0}\beta_{2,j}e_{2,j}$. We assume that $v$ is already of this form. Then
\begin{eqnarray*}
w&=&(b\sigma)^2(v)-tv\\
&=& \beta_{2,0}^{\sigma^2}e_{2,2}+ (\beta_{2,1}^{\sigma^2}-\beta_{2,0})e_{2,3}+\dotsm\in b\sigma(\Lambda)+t\Lambda.
\end{eqnarray*}
We repeat the modifications of $v$ similarly for $w$: By dividing $w$ by $\beta_{2,0}^{\sigma^2}$ and after subtracting a suitable linear combination of the $(b\sigma)^lw$, we obtain that $e_{2,2}\in b\sigma(\Lambda)+t\Lambda$. Hence also $e_{2,j}=(b\sigma)^{j-2}(e_{2,2})\in b\sigma(\Lambda)+t\Lambda$ for all $j\geq 2$. We can modify $v$ by a suitable linear combination of these vectors to obtain an element of the form $$v'=e_{1,0}+a_0e_{2,0}+a_1e_{2,1}\in \Lambda\setminus (b\sigma(\Lambda)+t\Lambda)$$ for some $a_0,a_1\in \overline{k}$. This implies $$(b\sigma)^2(v')-a_0^{\sigma^2}e_{2,2}-a_1^{\sigma^2}e_{2,3}=e_{1,2}\in \Lambda.$$ Finally we obtain
\begin{equation}\label{glMbesch}
\Lambda=\langle e_{1,0}+a_0e_{2,0}+a_1e_{2,1},e_{1,1}+a_0^{\sigma}e_{2,1},e_{i,j}\mid j\geq 2\rangle_{\mathcal{O}_L}
\end{equation}
for some $a_0,a_1\in\overline{k}$. We define a morphism $\mathbb{A}^2\rightarrow X_{\mu}(b)$ by mapping $(a_0,a_1)$ to the lattice in (\ref{glMbesch}). We choose $\{s_1^a\circ s_2^{a'}\mid (a,a')\in\mathbb{Z}^2\}$ as a set of representatives of $J/(K\cap J)\cong \mathbb{Z}^2$ in $J$. Here the $s_i\in GL_{i+1}\subset M_{\sb}= GL_2\times GL_3$ are as in the proof of Proposition \ref{propsurj} for superbasic $b$. As $X_{\mu}(b)$ is invariant under $J$, this yields a morphism $$\coprod_{\mathbb{Z}^2}\mathbb{A}^2\rightarrow X_{\mu}(b).$$ The morphism defines a bijection of geometric points which implies $$\pi_0(X_{\mu}(b))\cong J/(K\cap J)\cong \mathbb{Z}^2.$$

Note that there are also examples of $(G,\mu,b)$ such that the action of $J$ on $\pi_0(X_{\mu}(b))$ has more than one orbit. For example, one can use the methods of \cite{dimdlv} to show that for $G=GL_3$, $\mu=(3,1,0)$ and $b$ superbasic of slope $\frac{4}{3}$, the affine Deligne-Lusztig variety is isomorphic to $\coprod_{\mathbb{Z}\sqcup\mathbb{Z}}\mathbb{A}^2$. In this case, the action of $J$ on $\pi_0$ has two orbits, one for each copy of $\mathbb{Z}$ in the index set.

\section{Zero-dimensional affine Deligne-Lusztig varieties} 

\begin{kor}\label{lem3}
Let $G$ be simple and $\P=\M\N$ a proper standard parabolic subgroup. Let $b\in\M(L)$ be basic in $G$. Let $\mu$ be $G$-dominant with $\kappa_{\M}(b)=\mu$. Then $[b]=[t^\mu]$ is central.
\end{kor}
\begin{proof}
This is a special case of Theorem \ref{thmzshk}(1). Indeed, the assumption that $b$ is basic implies that $\langle\alpha_i,\nu\rangle=0$ for every simple root $\alpha_i$ of $A$ in $G$. Thus $(\mu,b)$ is indecomposable with respect to the Hodge-Newton decomposition.
\end{proof}

\begin{proof}[Proof of Theorem \ref{thmdim0}]
Note that $X_{\mu}(b)$ is open in $X_{\preceq\mu}(b)$. Thus an isolated point of the former variety is still isolated in the latter. But $J$ acts transitively on $\pi_0(X_{\preceq\mu}(b))$. Thus the first assertion implies the second. The inverse implication is trivial.

It remains to prove that the second and third assertion are equivalent. By Lemma \ref{lemnormgm}, $\dim X_{\mu}^G(b)=0$ if and only if all $X^{G_i}_{\mu_i}(b_i)$ are zero-dimensional. Here, the $G_i$ are the simple factors of $G_{\ad}$. Hence we may assume that $G$ is adjoint and simple. 

To prove that the second assertion implies the third assume that $\dim X_{\mu}(b)=0$. Let $x\in X_{\mu}(b)$. By Proposition \ref{propsurj} its connected component in $X_{\preceq\mu}(b)$ contains an element of $X_{\mu_{\min}}(b)$. But $X_{\mu}(b)$ is open in $X_{\preceq\mu}(b)$ and zero-dimensional. Thus $x\in X_{\mu_{\min}}(b)$ and $\mu=\mu_{\min}$.

In Section \ref{secnot} we saw that there is an $M_{\sb}$-dominant $\mu_{M_{\sb}}$ with $\mu_{M_{\sb},\dom}\preceq\mu$ and $X_{\mu_{M_{\sb}}}^{M_{\sb}}(b)\neq\emptyset$. If $\mu_{M_{\sb}}$ is not $G$-dominant, then there is a simple root $\alpha$ of $A$ in $N_{\sb}$ such that $\langle\alpha,\mu_{M_{\sb}}\rangle<0$. But then the proof of Proposition \ref{props0} provides a one-dimensional family in $X_{\mu}(b)$. (In this specific case, we may choose $\mu=\mu_{\min}$ in Proposition \ref{props0}.) This is a contradiction to $\dim X_{\mu}(b)=0$. Thus $\mu_{M_{\sb}}$ is $G$-dominant, and $\mu_{M_{\sb}}=\mu=\mu_{\min}$. As $\kappa_{M_{\sb}}(b)=\mu_{M_{\sb}}=\mu$, we have $\kappa_{\M}(b)=\mu$ in $\pi_1(\M)$ for every Levi subgroup $\M\supseteq M_{\sb}$.

We first apply this to $\M=\M_b$. As $(\mu,b)$ is indecomposable with respect to the Hodge-Newton decomposition, $\kappa_{\M_b}(b)=\mu$ implies that $\M_b=G$. Hence $b$ is basic in $G$. We consider two cases:

{\bf Case 1:} $b$ is superbasic. Then $G=PGL_h$ for some $h$ (compare \cite{GHKR}, 5.9). Write $\mu=(1,\dotsc,1,0,\dotsc,0)\in\mathbb{Z}^h/\mathbb{Z}$ with multiplicities $m$ and $h-m$. As $b$ is superbasic, $m$ and $h$ are coprime. The dimension formula for affine Deligne-Lusztig varieties reads in this special case
\begin{equation}\label{gldimdlv}
\dim X_{\mu}(b)=\langle\rho,\mu-\nu\rangle-\frac{h-1}{2}
\end{equation}
(compare \cite{GHKR}, \cite{dimdlv}). From this, one easily obtains $\dim X_{\mu}(b)=\frac{(m-1)(h-m-1)}{2}$. It vanishes if and only if $m$ or $h-m$ is equal to $1$.

{\bf Case 2:} There is a proper Levi subgroup $\M\subset G$ containing $b$. The observation made above implies that $\kappa_{\M}(b)=\mu$. Then the claim follows from Corollary \ref{lem3}.

The proof that the last assertion implies the second is immediate from Theorem \ref{thmzshk}(3) and the dimension formula for superbasic $b$, see Case 1 above. 
\end{proof}

\end{document}